\documentclass{article}      

\usepackage{epsfig}

\title{Golden Quantum Oscillator and Binet-Fibonacci Calculus}

\author{Oktay K. Pashaev and Sengul Nalci\\Department of Mathematics,\\ Izmir Institute of Technology \\ Urla-Izmir, 35430, Turkey}

\begin{document}
\newcommand{\be}{\begin{equation}}
\newcommand{\ee}{\end{equation}}
\newcommand{\bea}{\begin{eqnarray}}
\newcommand{\eea}{\end{eqnarray}}
\newcommand{\disp}{\displaystyle}
\newcommand{\la}{\langle}
\newcommand{\ra}{\rangle}

\newtheorem{thm}{Theorem}[subsection]
\newtheorem{cor}[thm]{Corollary}
\newtheorem{lem}[thm]{Lemma}
\newtheorem{prop}[thm]{Proposition}
\newtheorem{defn}[thm]{Definition}
\newtheorem{rem}[thm]{Remark}
\newtheorem{prf}[thm]{Proof}

\maketitle


\begin{abstract}
The Binet-Fibonacci formula for Fibonacci numbers is treated as a q-number(and q-operator) with Golden ratio bases $q=\varphi$ and $Q=-1/\varphi$. Quantum harmonic oscillator for this
Golden calculus is derived so that its spectrum is given just by Fibonacci numbers. Ratio of successive energy levels is found as the Golden sequence
and for asymptotic states it appears as the Golden ratio. This why we called this oscillator as the Golden oscillator. By double Golden bosons,
the Golden angular momentum and its representation in terms of Fibonacci numbers and the Golden ratio are derived.

\end{abstract}

\section{Introduction}
Fibonacci numbers are known from ancient time and have many applications from human proportions, architecture (Golden section), natural plants (branches of trees, arrangement of leaves) up to financial market \cite{Koshy}.

The numbers satisfy the recursion relation
$$F_1= F_2= 1\,\,\, {\rm (Initial\,\, Condition)},$$
$$F_n = F_{n-1}+F_{n-2}, \,\,\, {\rm for}\,\, n\geq 2 \,\,\, ({\rm Recursion\,\, Formula}).$$
First few Fibonacci numbers are $1,1,2,3,5,8,13,...$
For these numbers, starting from de Moivre, Lame and Binet, next representation is known as the Binet formula \cite{Koshy}:
\be  F_n = \frac{\varphi ^n- {\varphi'}^n}{\varphi- \varphi'},\ee where $\varphi, \, \varphi'$ are positive and negative roots of equation $$x^2-x-1=0.$$
These roots are explicitly \be  \varphi= \frac{1+\sqrt{5}}{2}, \,\,\,\,\,\,\,\,\, \varphi'= \frac{1-\sqrt{5}}{2}= -\frac{1}{\varphi}. \ee
Number $\varphi$ is known as the Golden ratio or the Golden section. There is a huge amount of works devoted to application of Golden ratio in many fields from natural phenomena to architecture and music.

Here we notice that Binet formula is a special realization of so-called $q$-numbers in $q$-calculus with two basis $q,\,\,Q,$ where $q= \varphi, \,\, Q= -\frac{1}{\varphi}= \varphi'.$ The $(Q,q)$ calculus generalizes the $q$-calculus. In particular cases when $Q=1$ it becomes non-symmetrical calculus. In case $Q= \frac{1}{q}$ it becomes so-called symmetrical $q$-calculus. It appears in study of generalized quantum $q$-harmonic oscillator \cite{Arik}, \cite{Chakrabarti} and mentioned as a convenient form for generalization, in book \cite{Kac}.

Recently, we found that it appears naturally in construction of $q$-Binomial formula for noncommutative elements. Noncommutative $q$-binomials where considered in \cite{Nalci} for description $q$-Hermite polynomial solutions for $q$-Heat equation. From another side it appears also in description of AKNS Hierarchy of integrable systems where $Q=R$ is recursion operator of AKNS Hierarchy and $q$ is the spectral parameter \cite{Pashaev}.

In the present article we like to explore possibility to interpret Binet formula for Fibonacci numbers as $q$-number with Golden section base $\varphi.$ Then Fibonacci numbers appear as $q$-numbers. We construct quantum harmonic oscillator for Golden $q$-Calculus case so that its spectrum is given by just Fibonacci numbers and
ratio of successive energy levels is given as the Golden sequence
and for asymptotic states it appears as the Golden ratio.
 Though the results for harmonic oscillator with generic $Q-q$ and its reductions to symmetrical and non-symmetrical cases are known
\cite{Arik1}, \cite{Biedenharn}, \cite{Macfarlane}, we think that the special case with Golden ratio base was not described before in literature. Due to importance and wide applicability of Fibonacci numbers in different fields, we think that explicit realization of it in form of quantum oscillator with Golden ratio base, which we call Golden quantum oscillator, could be of interest. Specially, a realization of this type of calculus could describe the Golden ratio in noncommutative geometry.

Finally, our Golden oscillator should not be mixed with the  Fibonacci Oscillator of \cite{Arik}
with generic bases $q_1$ and $q_2$, thought it is a particular case of it. In that papers no relations with Golden ratio and Binet formula as well as asymptotic properties of energy levels are discussed.
\section{Golden $q$-Calculus}
In $(Q,q)$ calculus we have number \be [n]_{Q,q}= \frac{Q^n-q^n}{Q-q}.\label{Qqnumber}\ee
If we choose $Q=\varphi=\frac{1+\sqrt{5}}{2}$ and $q=\varphi'=\frac{1-\sqrt{5}}{2}=-\frac{1}{\varphi}.$
Then (\ref{Qqnumber}) becomes Binet's formula for Fibonacci numbers as $(\varphi,\, \varphi')$-numbers :
 \be  F_n = \frac{\varphi ^n- {\varphi'}^n}{\varphi- \varphi'}= [n]_{\varphi,\varphi'}\equiv [n]_{F}.\label{binet}\ee
This definition can be extended to arbitrary real number $x$,
\be [x]_{\varphi,\varphi'}\equiv [x]_F=  \frac{\varphi ^x- {\varphi'}^x}{\varphi- \varphi'}= \frac{\varphi ^x- (-\frac{1}{\varphi})^x}{\varphi+ \frac{1}{\varphi}}\equiv F_x, \ee though due to negative sign for the second base,  it is not a real number for general $x$
 \be F_x =\frac{1}{\varphi +
  \frac{1}{\varphi}} \left(\varphi^x - e^{i\pi x} \frac{1}{\varphi^x}\right)=\frac{1}{\sqrt{5}} \left(\varphi^x - e^{i\pi x} \frac{1}{\varphi^x}\right) .\ee
 Instead of  real number $x$ we can also consider complex numbers $z= x+i y,$

\textbf{Example :} it is easy to see that
$$\lim_{n\rightarrow \infty } \frac{[n+1]_F}{[n]_F}=\lim_{n\rightarrow \infty } \frac{F_{n+1}}{F_n}=\varphi.$$
The addition formula for Golden numbers is given in the form
\be [n+m]_F= F_{n+m}={\varphi}^n F_m + \left({-\frac{1}{\varphi}}\right)^m F_n \label{multipleaddition}. \ee
By using (\ref{binet}) we can get
\be \varphi^N= \varphi F_N+ F_{N-1}, \,\,\,\,\,\,\,\,\,\,\,\, {\varphi'}^N= \varphi' F_N+ F_{N-1}, \label{Fibonacciprop}\ee and the above formula (\ref{multipleaddition}) can be rewritten as
\bea F_{n+m}&=& F_n F_{m-1}+ F_{n+1} F_m \nonumber \\
&=& F_{n-1} F_m+ F_n F_{m+1}.\eea

The substraction formula can be obtained from it by changing $m\rightarrow -m$ as
 \be F_{n-m}=[n-m]_{F}= {\varphi}^{n} [-m]_{F}+\left({-\frac{1}{\varphi}}\right)^{-m} [n]_{F} \ee
 or by using the equality
 $$[-n]_{F}=-(-1)^{-n} [n]_{F}$$ it can be also written
\bea [n-m]_{F}&=&\left({-\frac{1}{\varphi}}\right)^{-m} ([n]_{F}-{\varphi}^{n-m} [m]_{F}) \nonumber \\
&=& \left({-\frac{1}{\varphi}}\right)^{-m} F_n-\frac{\varphi^n}{(-1)^m} F_m. \eea
or
\bea F_{n-m}
&=& \left({-\frac{1}{\varphi}}\right)^{-m} F_n-\frac{\varphi^n}{(-1)^m} F_m. \eea
\textbf{DEFINITION (Higher Fibonacci Numbers):}
\be F_n^{(m)}\equiv \frac{(\varphi^m)^n-(\varphi'^m)^n}{\varphi^m-\varphi'^m}= [n]_{\varphi^m,\varphi'^m}\ee and $F_n^{(1)}\equiv F_n.$

By definition, the multiplication rule for Golden numbers is given by next formula
\be[n m]_{\varphi,-\frac{1}{\varphi}}=F_{n m}= [n]_{\varphi,-\frac{1}{\varphi}} [m]_{\varphi^n,{(-\frac{1}{\varphi})}^n}=F_n F_m^{(n)},\ee
and the division rule is
\bea   \left[\frac{m}{n}\right]_{\varphi,\varphi'}&=&\frac{[m]_{\varphi,\varphi'}}{[n]_{\varphi^{m/n},\varphi'^{m/n}}}= \frac{[m]_{\varphi^{1/n},\varphi'^{1/n}}}{[n]_{\varphi^{1/n},\varphi'^{1/n}}}\nonumber \\
F_{\frac{m}{n}}&=& \frac{F_m}{F_n^{(\frac{m}{n})}}= \frac{F_m^{(\frac{1}{n})}}{F_n^{(\frac{1}{n})}}.  \eea

Higher Fibonacci numbers can be written in terms of ratio of Fibonacci numbers as follows
\be F_n^{(m)}= \frac{F_{m n}}{F_m}.\ee
From definition (\ref{binet}) we have the following relation
\be F_{-n}= (-1)^{n+1} F_n.\ee

For any real $x,y$ \bea [x+y]_F &=& \varphi^x [y]_F+\left(-\frac{1}{\varphi}\right)^y [x]_F \nonumber \\
&=& \varphi^y [x]_F+\left(-\frac{1}{\varphi} \right)^x [y]_F \eea
which are written in terms of Fibonacci numbers as follows
\bea F_{x+y}&=&\varphi^x F_y + (-\frac{1}{\varphi})^y F_x \nonumber \\
&=& \varphi^y F_x+\left(-\frac{1}{\varphi} \right)^x F_y . \eea
For real $x,$ we have the Fibonacci recurrence relation
\be [x]_F=[x-1]_F+[x-2]_F \Rightarrow F_x=F_{x-1}+F_{x-2}. \ee

\textbf{Example :} Golden $\pi$
$$F_\pi= [\pi]_F \simeq= 4,73068+0,0939706i$$
\section{Golden Derivative}
The Fibonacci or Golden derivative we define as an operator
\be  F_{x \frac{d}{dx}} = \frac{\varphi^{x \frac{d}{dx}}- \varphi'^{x \frac{d}{dx}}}{\varphi - \varphi'} = [x \frac{d}{dx}]_F . \ee
Then Golden derivative for any function $f(x)$ is given as
 \be F_{x \frac{d}{dx}}f(x) = D_F f(x)= \frac{f(\varphi x)-f(-\frac{x}{\varphi})}{(\varphi+\frac{1}{\varphi})x}= \frac{\left( M_\varphi-M_{-\frac{1}{\varphi}}\right) f(x)}{(\varphi+\frac{1}{\varphi}x)} \label{goldenderivative} .\ee
 Here, arguments are scaled by Golden ratio: $x\rightarrow \varphi x$ and $x\rightarrow -\frac{x}{\varphi}$. It can be written in terms of Golden ratio dilatation operator  $$M_{\varphi} f(x)= f(\varphi x),$$
 where  $f(x)$- smooth function and its operator form can also be written as $$M_\varphi= \varphi^{x \frac{d}{d x}}= \left(\frac{1+\sqrt{5}}{2} \right)^{x \frac{d}{d x}}. $$

Function $A(x)$ we called the Golden periodic function if
\be D_F A(x) = 0 \ee
which implies
\be A(\varphi x) = A(- \frac{1}{\varphi} x) .\ee
As an example we have:
\be A(x) = \sin \left(\frac{\pi}{\ln \varphi} \ln |x|\right)\ee

\textbf{Example \,\,1:} Application of Golden derivative operator $D_F$ on $x^n$ gives $$D_F x^n= F_n x^{n-1}$$ or $$F_n= \frac{D_F x^n}{x^{n-1}},$$
so it generates Fibonacci numbers.

\textbf{Example \,\,2:} $$D_F e^x = \sum_{n=0}^\infty \frac{F_n}{n!} x^n$$ or $$D_F e^x= \frac{e^{\varphi x}-e^{-\frac{x}{\varphi}}}{\varphi+\frac{1}{\varphi}}= \frac{2 e^{\frac{x}{2}} \sinh {\frac{\sqrt{5}}{2}x}}{\sqrt{5} x}=\sum_{n=0}^\infty \frac{F_n}{n!} x^n .$$
For $x=1$ it gives next summation formula
\be \sum_{n=0}^\infty \frac{F_n}{n!} =  e^{\frac{1}{2}}\frac{ \sinh {\frac{\sqrt{5}}{2}}}{\frac{\sqrt{5}}{2} } . \ee
\subsection{Golden Leibnitz Rule}
 We derive the Golden Leibnitz rule
 \be D_F (f(x) g(x))= D_F f(x) g(\varphi x)+f(-\frac{x}{\varphi}) D_F g(x) \label{multilebnitz}.\ee
 By symmetry the second form of the Leibnitz rule can be derived as
 \be D_F (f(x) g(x))= D_F f(x) g(-\frac{x}{\varphi})+f(\varphi x) D_F g(x) \label{multilebnitz2}.\ee
 These formulas can be  rewritten  in explicitly symmetrical form :
\be D_F (f(x) g(x))=D_F f(x) \left(\frac{g(\varphi x)+g(-\frac{x}{\varphi})}{2} \right) + D_F g(x) \left(\frac{f(\varphi x) + f(-\frac{x}{\varphi})}{2} \right). \ee
More general form of Golden Leibnitz formula is given with arbitrary $\alpha$ ,
$$D_F (f(x) g(x))= \left( \alpha f(-\frac{x}{\varphi}))+(1-\alpha) f(\varphi x)\right) D_F g(x)+ \left(\alpha g(\varphi x)+(1-\alpha)g(-\frac{x}{\varphi}))\right)D_F f(x)$$

Now we may compute the golden derivative of the quotient of $f(x)$ and $g(x)$. From (\ref{multilebnitz}) we have
\be D_F \left( \frac{f(x)}{g(x)} \right)= \frac{D_F f(x) g(\varphi x)-D_F g(x) f(\varphi x)}{g(\varphi x) g( -\frac{x}{\varphi})} \label{1} .\ee
However, if we use (\ref{multilebnitz2}), we get
\be D_F \left( \frac{f(x)}{g(x)} \right)= \frac{D_F f(x) g(-\frac{x}{\varphi})-D_F g(x) f( -\frac{x}{\varphi})}{g(\varphi x) g(-\frac{x}{\varphi})}. \label{2} \ee

In addition to the formulas (\ref{1}) and (\ref{2}) one may determine one more representation in symmetrical form
\be D_F \left( \frac{f(x)}{g(x)} \right)= \frac{1}{2}\frac{D_F f(x) (g(-\frac{x}{\varphi})+ g(\varphi x))-D_F g(x) (f( -\frac{x}{\varphi}) + f(\varphi x))}{g(\varphi x) g(-\frac{x}{\varphi})} . \label{12} \ee
In particular applications one of these forms could be more useful than others.
\subsection{Golden Taylor Expansion}
\begin{thm} Let the Golden derivative operator $D_F$ is a linear operator on the space of polynomials, and $$P_n (x)\equiv \frac{x^n}{F_n!}\equiv \frac{x^n}{F_1 F_2 ...F_n}$$ satisfy the following conditions : \\
(i)\,\, $P_0(0)= 1$ and $P_n(0)=0$ for any  $n\geq 1$ ; \\
(ii) deg $P_n=n$ ; \\
(iii) $D_F P_n (x)= P_{n-1} (x)$ for any  $n\geq 1 $ , and $D_F(1)=0$
Then, for any polynomial $f(x)$ of degree N,one has the following Taylor formula :
$$ f(x)=\sum_{n=0}^N (D_F^n f) (0) P_n (x)= \sum_{n=0}^N (D_F^n f) (0) \frac{x^n}{F_n!}.$$
\end{thm}
In the limit $N\rightarrow \infty$ (when it exists) this formula can determine some new function \be f_F(x)= \sum_{n=0}^\infty (D_F^n f)(0) \frac{x^n}{F_n!}\ee
which we can call the Golden (or Fibonacci) function.

\textbf{Example :\, (Golden Exponential) }
The Golden exponential functions are
\be e_F^x\equiv \sum_{n=0}^\infty \frac{x^n}{F_n!}\,\, ; \,\,\,\,\,\, E_F^x \equiv \sum_{n=0}^\infty (-1)^{\frac{n(n-1)}{2}} \frac{x^n}{F_n!},\ee
and for $x=1,$ we get the Fibonacci natural base as follows
$$e_F^x\equiv \sum_{n=0}^\infty \frac{1}{F_n!} \equiv e_F.$$
These functions are entire analytic functions. For the second function explicitly we have
\be E_F^x = 1 + \frac{x}{F_1 !} - \frac{x^2}{F_2 !} - \frac{x^3}{F_3 !} + \frac{x^4}{F_4 !} + \frac{x^5}{F_5 !} - \frac{x^6}{F_6 !} - \frac{x^7}{F_7 !}
+ \frac{x^8}{F_8 !} + \frac{x^9}{F_9 !} -...\ee
The Golden derivative of these exponential functions are found
$$D_F e_F^{k x}= k e_F^{k x},$$
$$D_F E_F^{k x}= k E_F^{-k x}$$
for arbitrary constant $k$ (or F-periodic function).
Then these two functions give the general solution of the hyperbolic F-oscillator
\be (D^2_F - k^2) \phi(x) = 0, \ee
as
\be \phi (x) = A e_F^{k x} + B e_F^{-k x}, \ee
and elliptic F-oscillator
\be (D^2_F + k^2) \phi(x) = 0, \ee
\be \phi (x) = A E_F^{k x} + B E_F^{-k x}.\ee
We have next Euler formulas
\be e_F^{ix} = \cos_F x + i \sin_F x ,\ee
\be E_F^{ix} = \cosh_F x + i \sinh_F x ,\ee
and relations
\be Cosh_F x = \cos_F x ,\ee
\be Sinh_F x = \sin_F x ,\ee
where
\be Cosh_F x \equiv \frac{E_F^x + E_F^{-x}}{2},\,\,\,\, Sinh_F x \equiv \frac{E_F^x - E_F^{-x}}{2} .\ee
We notice here that these relations are valid due to alternating character of second exponential function.

\textbf{Example :\, (F-Oscillator) }

For F-oscillator
\be D^2_F x + \omega^2 x = 0\ee
the general solution is
\be x(t) = a E^{\omega t}_F + b E^{-\omega t}_F = a' Cosh_F \omega t + b' Sinh_F \omega t = a' \cos_F \omega t + b' \sin_F \omega t \ee

\subsection{Golden Binomial}
Golden Binomial we define  as
\be (x+y)_F^n = (x+\varphi^{n-1} y)(x-\varphi^{n-3} y)...(x+ (-1)^{n-1}\varphi^{-n+1} y)\ee and it has n-zeros at  the Golden ratio powers
$$\frac{x}{y}=-\varphi^{n-1},\,\,\,\, \frac{x}{y}=-\varphi^{n-3},\,\,\,\,...,\frac{x}{y}=-\varphi^{-n+1}.$$

For Golden binomial next expansion is valid
\bea (x+y)_F^n \equiv (x+y)_{\varphi,-\frac{1}{\varphi}}^n &=&\sum^{n}_{k=0}{ n \brack k}_{F} (-1)^{\frac{k(k-1)}{2}} x^{n-k} y^k  \nonumber \\
&=& \sum^{n}_{k=0} \frac{F_n!}{F_{n-k}! F_k!}(-1)^{\frac{k(k-1)}{2}}  x^{n-k} y^k \label{goldenbinomexpansion}\eea
The proof is easy by induction.

Application of Golden derivative to the Golden binomial gives
$$D_F^x (x+y)_F^n=F_n (x+y)_F^{n-1},$$
$$D_F^y (x+y)_F^n=F_n (x-y)_F^{n-1} .$$
It means
$$D_F^x \frac{(x+y)_F^n}{F_n!}=\frac{(x+y)_F^{n-1}}{F_{n-1}!} ,$$
$$D_F^y \frac{(x+y)_F^n}{F_n!}=\frac{(x-y)_F^{n-1}}{F_{n-1}!} ,$$

$$\left( D_F^y\right)^{n} (x+y)_F^{n} .$$
For $n=2k$ we have
$$\left( D_F^y\right)^{2 k} (x+y)_F^{2 k}= (-1)^k F_{2 k}!,$$
and for $n=2k+1$ we get  $$\left( D_F^y\right)^{2k+1} (x+y)_F^{2k+1}= (-1)^k F_{2k+1}!$$

In terms of Golden binomial we introduce the Golden polynomials

\be P_n (x) = \frac{(x-a)_F^n}{F_n!}\ee
where $n=1,2,...$, and $P_0(x) =1$ with property
 \be D_F^x P_n(x) = P_{n-1}(x) .\ee
For even and odd polynomials we have products
\be P_{2n} (x) = \frac{1}{F_{2n}!} \prod^n_{k=1} (x- (-1)^{n+k}\varphi^{2k-1} a) (x + (-1)^{n+k}\varphi^{-2k +1} a) ,\ee
\be P_{2n+1} (x) = \frac{(x - (-1)^n a)}{F_{2n+1}!} \prod^n_{k=1} (x- (-1)^{n+k}\varphi^{2k} a) (x - (-1)^{n+k}\varphi^{-2k} a) .\ee

By using (\ref{Fibonacciprop}) it is easy to find
\be \varphi^{2k} + \frac{1}{\varphi^{2k}} = F_{2k} + 2 F_{2k-1} ,\ee
\be \varphi^{2k+1} - \frac{1}{\varphi^{2k+1}} = F_{2k+1} + 2 F_{2k} .\ee
Then we can rewrite our polynomials in terms of just Fibonacci numbers
\be P_{2n} (x) = \frac{1}{F_{2n}!} \prod^n_{k=1} (x^2 - (-1)^{n+k} (F_{2k-1} + 2 F_{2k-2})x a - a^2) ,\ee
\be P_{2n+1} (x) = \frac{(x - (-1)^n a)}{F_{2n+1}!} \prod^n_{k=1} (x^2- (-1)^{n+k}(F_{2k} + 2 F_{2k-1})x a + a^2) .\ee
First few polynomials are
\be P_1(x) = (x-a)\ee
\be P_3 (x) = \frac{1}{2} (x+a)(x^2 - 3 x a + a^2 )\ee
\be P_5 (x) = \frac{1}{2\cdot 3 \cdot 5} (x-a)(x^2 + 3 x a + a^2 )(x^2 - 7 x a + a^2 )\ee
\be P_7 (x) = \frac{1}{2\cdot 3 \cdot 5 \cdot 8 \cdot 13} (x+a)(x^2 - 3 x a + a^2 )(x^2 + 7 x a + a^2 )(x^2 - 18 x a + a^2 )\ee

\be ... \ee
\be P_2 (x) =  (x^2 -  x a - a^2 )\ee
\be P_4 (x) = \frac{1}{2\cdot 3} (x^2 +  x a - a^2 )(x^2 - 4 x a - a^2)\ee
\be P_6 (x) = \frac{1}{2\cdot 3\cdot 5 \cdot 8} (x^2 -  x a - a^2 )(x^2 + 4 x a - a^2)(x^2 - 11 x a - a^2)\ee
\be ... \ee
\subsection{Noncommutative Golden Ratio and Golden Binomials}

By choosing $q=-\frac{1}{\varphi}$ and $Q=\varphi,$ in general Q-commutative q-binomial \cite{Nalci1}, where $\varphi$ is the Golden section, we obtain the Binet-Fibonacci Binomial formula with Golden non-commutative plane $(y x= \varphi x y).$ (It should be compared with Golden ratio $b=\varphi a$).
\bea (x+y)^n_{-\frac{1}{\varphi}} &=& (x+y) (x+(-\frac{1}{\varphi})y) (x+(-\frac{1}{\varphi})^2 y)...(x+(-\frac{1}{\varphi})^{n-1}y) \nonumber \\
&=&\sum_{k=0}^n {n \brack k }_{\varphi,-\frac{1}{\varphi}} (-\frac{1}{\varphi})^{\frac{k(k-1)}{2}} x^{n-k} y^k \nonumber \\
&=& \sum_{k=0}^n \frac{F_n !}{F_k ! F_{n-k}!} \left(-\frac{1}{\varphi}\right)^{\frac{k(k-1)}{2}} x^{n-k} y^k ,\eea
where $F_n$ are Fibonacci numbers.

\subsection{Golden Pascal Triangle}

The Golden binomial coefficients are defined by
\be { n \brack k}_F= \frac{[n]_{F}!}{[n-k]_{F}! [k]_{F}!}= \frac{F_n!}{F_{n-k}! F_k!} \label{goldenbinom}\ee
with $n$ and $k$ being nonnegative integers, $n\geq k$ and are called the Fibonomials.
Using the addition formula for Golden numbers (\ref{multipleaddition}), we write following expression
$$F_n=F_{n-k+k}=(-\frac{1}{\varphi})^k F_{n-k}+\varphi^{n-k} F_k,$$
and from (\ref{Fibonacciprop}) it can be written as follows
\bea F_n &=& F_{n-k-1} F_k + F_{n-k} F_{k+1} \nonumber \\
&=& F_{n-k} F_{k-1}+ F_{n-k+1} F_k. \eea

With the above definition (\ref{goldenbinom})we have next recursion formulas
\bea { n \brack k}_{F}&=& \frac{(-\frac{1}{\varphi})^k [n-1]_{F}!}{[k]_{F}! [n-k-1]_{F}!} + \frac{\varphi^{n-k} [n-1]_{F}!}{[n-k]_{F}![k-1]_{F}!} \nonumber \\
&=& (-\frac{1}{\varphi})^k { n-1 \brack k}_{F} + \varphi^{n-k} { n-1 \brack k-1}_{F} \label{goldenpascal1}\\
&=& \varphi^k { n-1 \brack k}_{F} + (-\frac{1}{\varphi})^{n-k} { n-1 \brack k-1}_{F} .\label{goldenpascal2} \eea
These two rules determine the multiple Golden Pascal triangle, where $1\leq k\leq n-1.$
Then, we can construct Golden Pascal triangle as follows

\[
\matrix{&&&&& 1 &&& &&\cr\cr &&&&\swarrow & &\searrow &\cr\cr
& &&1& &&& 1& &\cr &&\swarrow&  &\searrow  -\frac{1}{\varphi} &&
\varphi \swarrow && \searrow&\cr\cr &1& &&& [2]_{F}& &&&
1&\cr \swarrow & & \searrow (-\frac{1}{\varphi})^2 &
& \varphi \swarrow & & \searrow  -\frac{1}{\varphi} & & \varphi^2
\swarrow & & \searrow \cr &\cdots&&&&\cdots&&&&\cdots&\cr}
\]

\subsection{Remarkable Limit}

From Golden binomial expansion (\ref{goldenbinomexpansion}) we have

\bea (1+y)_F^n  &=&\sum^{n}_{k=0}{ n \brack k}_{F} (-1)^{\frac{k(k-1)}{2}}  y^k  \nonumber \\
&=& \sum^{n}_{k=0} \frac{F_n!}{F_{n-k}! F_k!}(-1)^{\frac{k(k-1)}{2}}   y^k  .\label{goldenbinomexpansion1}\eea
Then
\be \left(1+ \frac{y}{\varphi^n}\right)_F^n
= \sum^{n}_{k=0} \frac{F_n!}{F_{n-k}! F_k!}(-1)^{\frac{k(k-1)}{2}}  \frac{y^k}{\varphi^{n k}} \label{goldenbinomexpansion2}\ee
or by opening Fibonomials and taking limit
\be \lim_{n \rightarrow \infty}\left(1+ \frac{y}{\varphi^n}\right)_F^n
= \sum^{\infty}_{k=0} \frac{1}{ F_k!}(-1)^{\frac{k(k-1)}{2}}  \frac{y^k}{\varphi^{\frac{k(k-1)}{2}} (\varphi + \frac{1}{\varphi})^k} \label{limit}\ee
\be \lim_{n \rightarrow \infty}\left(1+ \frac{y}{\varphi^n}\right)_F^n
= \sum^{\infty}_{k=0} \frac{1}{ [k]_{-\varphi^2}!}
 \left( \frac{y \varphi}{\varphi^2 +1} \right)^k \label{limit1}\ee
where we introduced $q$-number, $[k]_q = 1+ q + ... + q^{k-1}$,  with base $q = - \varphi^2$, so that
\be [k]_{-\varphi^2} = 1+ (-\varphi^2) + ... + (-\varphi^2)^{k-1} = \frac{(-\varphi^2)^k - 1}{(-\varphi^2) - 1} .\ee
The last expression allow us to rewrite the limit in terms of Jackson q-exponential function $e_q(x)$ with $q = - \varphi^2$,
\be \lim_{n \rightarrow \infty}\left(1+ \frac{y}{\varphi^n}\right)_F^n
= e_{- \varphi^2}\left( \frac{y \varphi}{\varphi^2 +1} \right) \label{limit2}\ee
or finally we have remarkable limit
\be \lim_{n \rightarrow \infty}\left(1+ \frac{y}{\varphi^n}\right)_F^n
= e_{- \varphi^2}\left( \frac{y }{\sqrt{5}} \right) \label{limit3} .\ee
In particular case it gives
\be \lim_{n \rightarrow \infty}\left(1+ \frac{\sqrt{5}}{\varphi^n}\right)_F^n
= e_{- \varphi^2} ( 1) .\label{limit4}\ee
\subsection{Golden Integral}
\subsubsection{Golden Antiderivative}
\begin{defn} The function $G(x)$ is Golden antiderivative of $g(x)$ if $D_{F} G(x)=g(x).$

It is denoted by
\be G(x)=\int g(x)d_{F} x.\ee

\end{defn}

$$D_{F} G(x)=0 \Rightarrow G(x)=C-constant$$
or
$$D_{F} G(x)=0 \Rightarrow G(\varphi x)=G( -\frac{x}{\varphi})$$ is called the Golden 'periodic' function.

\subsubsection{Golden-Jackson Integral}
By inverting equation (\ref{goldenderivative}) and expanding inverse operator we find Jackson type representation for anti-derivative.

\be G(x)=\int g \left(\frac{x}{\varphi}\right) d_{Q} x =(1-Q) x \sum_{k=0}^\infty Q^k f \left(\frac{x}{\varphi} Q^k \right)\ee where $Q\equiv -\frac{1}{\varphi^2}.$

\section{Golden oscillator }
Now we construct quantum oscillator with spectrum in the form of Fibonacci numbers. Since in this oscillator the base in commutation relations is $\varphi$-Golden ratio, we called it as Golden oscillator.
The algebraic relations for Golden Oscillator are
\be b b^+ - \varphi b^+ b= (-\frac{1}{\varphi})^N \label{alrel1}\ee or
\be b b^+ + \frac{1}{\varphi} b^+ b= \varphi^N,\label{alrel2} \ee
where $N$ is the hermitian number operator and $\varphi$ is the deformation parameter. The bosonic Golden-oscillator is defined by three operators $b^{+}$, $b$ and $N$ which satisfy the  commutation relations:
\be [N,b^{+}]=b^{+}, \,\,\,\,\,\,\,\,\,\,\, [N,b]=-b. \label{multiplecomrel}\ee
By using the definition of number operator with basis $\varphi$ we find following equalities
\be [N+1]_F -\varphi [N]_F=(-\frac{1}{\varphi})^N\ee
\be[N+1]_F + \frac{1}{\varphi} [N]_F=\varphi^N ,\ee
where $$[N]_F=\frac{\varphi^N-(-\frac{1}{\varphi})^N}{\varphi+\frac{1}{\varphi}}$$
is the Fibonacci number operator. Here operator $(-1)^N = e^{i \pi N}$.

By comparison the above operator relations with algebraic relations (\ref{alrel1}) and (\ref{alrel2}) we have
$$b^{+} b=[N]_F, \,\,\,\,\,\,\,\,b b^{+}=[N+1]_F .$$
Here we should note that the number operator $N$ is not equal to $b^{+} b$ as in ordinary case.
By using the property of Fibonacci numbers (\ref{Fibonacciprop}) the algebraic relations (\ref{alrel1}) and (\ref{alrel2}) are equivalent to Fibonacci rule for operators $$F_{N+1}= F_N+ F_{N-1}.$$

\begin{prop}\bea \left[[N]_F,b^{+}\right]&=& \{[N]_F-[N-1]_F \}b^{+} \nonumber \\
&=& b^{+}\{[N+1]_F-[N]_F \} \label{multiple6}\eea
\end{prop}
\begin{prop} we have following equality for $n=0,1,2,..$
\be [[N]_F^n,b^{+}]=\{[N]_F^n-[N-1]_F^n\} b^{+}\ee
\end{prop}
\begin{prf}By using mathematical induction to show the above equality is not difficult.\\
\end{prf}
\begin{cor}For any  function expandable to power series (analytic) $F(x)=\sum_{n=0}^\infty c_n x^n$ we have the following relation
\bea [F([N]_F),b^{+}]&=&\{F([N]_F)-F([N-1]_F)\}b^{+} \nonumber \\
&=&b^{+}\{F([N+1]_F)-F([N]_F)\} \eea
and
\be b^+ F([N+1]_F)= F([N]_F) b^+ \label{a}\ee
or
\be F(N)b^{+}=b^{+}F(N+1) .\label{relation1} \ee
\end{cor}

By using the eigenvalues of the Number operator $$N|n\ra_F=n|n\ra_F,$$
$$[N]_F|n\ra_F=F_N|n\ra_F=
[n]_F |n\ra_F=F_n|n\ra_F$$
we get Fibonacci numbers as eigenvalues of $[N]$-operator, where we call $F_N$ as Fibonacci operator and we denote $|n\ra_{\varphi,-\frac{1}{\varphi}}\equiv |n\ra_F.$\\
The basis of the Fock space is defined by repeated action of the creation operator $b^+$ on the vacuum state, which is annihilated by $b|0\ra_F=0$
\be |n \ra_F= \frac{(b^{+})^n}{\sqrt{F_1 \cdot F_2\cdot ...F_n}}|0 \ra_F,\ee
where $[n]_F! = F_1 \cdot F_2\cdot ...F_n.$

In the limit $$\lim_{n \rightarrow \infty} \frac{F(n+1)}{F(n)}=\lim_{n \rightarrow \infty} \frac{[n+1]_F}{[n]_F}=\frac{1+\sqrt{5}}{2}\equiv \varphi=\approx 1,6180339887,$$ which is the Golden ratio.

The number operator $N$ for Fibonacci case is written in two different forms according to even or odd eigenstates $N|n\ra_F= n|n\ra_F.$
For $n=2k,$ we get
\be N= \log_{\varphi} \left(\frac{\sqrt{5}}{2} F_N+\sqrt{\frac{5}{4} F_N^2 +1} \right),\ee
and for $n=2k+1,$
\be N= \log_{\varphi} \left(\frac{\sqrt{5}}{2} F_N-\sqrt{\frac{5}{4} F_N^2 -1} \right),\ee
where $[N]_F$ is Fibonacci number operator defined as
$$[N]_F= \frac{\varphi^N-(-\frac{1}{\varphi})}{\varphi-(-\frac{1}{\varphi})}=F_N.$$

As a result, the Fibonacci numbers are the example of $(q,Q)$ numbers with two basis and one of the base is Golden Ratio.This is why we called the corresponding $q-$ oscillator as a Golden oscillator or Binet-Fibonacci Oscillator.
The Hamiltonian for $q$-Binet-Fibonacci oscillator is written as a Fibonacci number operator $$H=\frac{\hbar \omega}{2} (b^+ b+bb^+)=\frac{\hbar \omega}{2}\left( [N+1]_F +[N]_F\right)=\frac{\hbar \omega}{2}F_{N+2}, $$ where $b b^+= [N+1]_F=F_{N+1}, \,\,\,\,\, b^+ b= [N]_F=F_{N}.$
According to the Hamiltonian, the energy spectrum of this oscillator is written in terms of Fibonacci numbers sequence,
$$E_n=\frac{\hbar \omega}{2} \left([n]_{\varphi,-\frac{1}{\varphi}}+[n+1]_{\varphi,-\frac{1}{\varphi}}\right)=\frac{\hbar \omega}{2}\left(F_n+F_{n+1}\right)=\frac{\hbar \omega}{2}F_{n+2},$$
$$E_n=\frac{\hbar \omega}{2}F_{n+2}.$$
A first energy eigenvalues
$$E_0= \frac{\hbar \omega}{2}F_2=\frac{\hbar \omega}{2}, $$ which is exactly the same ground state as in the ordinary case. Higher energy excited states are given by Fibonacci sequence
$$E_1=\frac{\hbar \omega}{2}F_3=\hbar \omega,\,\,\,\,\,\,\,E_2= \frac{3 \hbar \omega}{2},\,\,\,\,\,\,\, E_3= \frac{5 \hbar \omega}{2},... $$
In Figure 1 we show the quantum Fibonacci tree for this oscillator.

\begin{figure}[h]
\begin{center}
\epsfig{figure=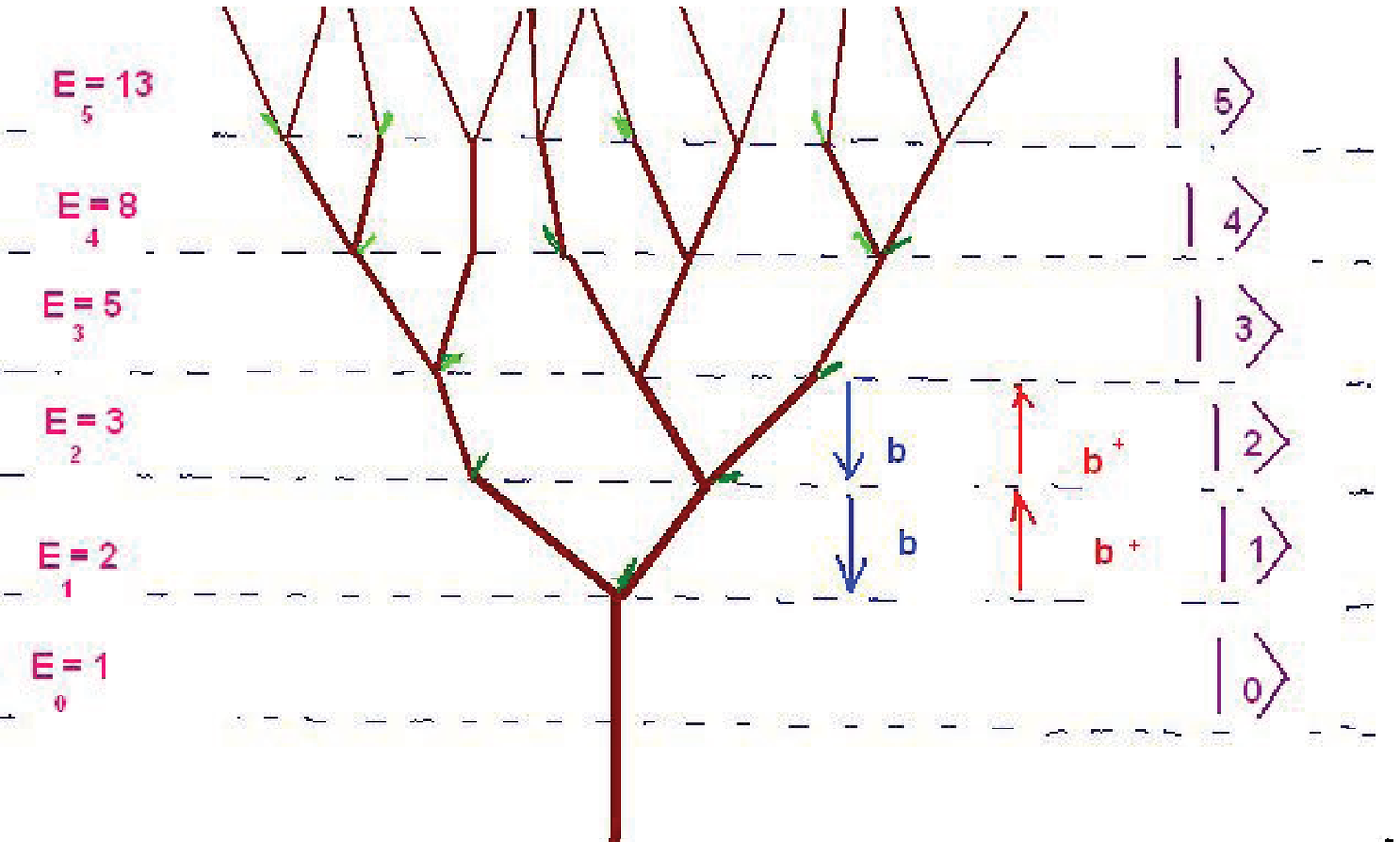,height=8cm,width=10cm}
\end{center}
\caption{Quantum Fibonacci Tree for Golden Oscillator}
\end{figure}

The difference between two consecutive energy levels of our oscillator is found as
$$\triangle E_n= E_{n+1}-E_n= \frac{\hbar \omega}{2} F_{n+1}.$$
Then the ratio of two successive energy levels $\frac{E_{n+1}}{E_n}$ gives the Golden sequence, and for the limiting case of higher excited states $n\rightarrow \infty$  it is the Golden ratio
$$\lim_{n \rightarrow \infty} \frac{E_{n+1}}{E_n}= \lim_{n \rightarrow \infty} \frac{F_{n+3}}{F_{n+2}}=\lim_{n \rightarrow \infty}\frac{[n+3]_F}{[n+2]_{F}}=\frac{1+\sqrt{5}}{2}= \varphi \approx 1,6180339887.$$
This property of asymptotic states to relate each other by a Golden ratio, leads us to call this oscillator as a Golden oscillator.

We have the following relations between $q$- creation and annihilation operators and standard creation and annihilation operators
\be b^+= a^+ \sqrt{\frac{F_{N+1}}{N+1}}= \sqrt{\frac{F_N}{N}} a^+\ee
\be b= \sqrt{\frac{F_{N+1}}{N+1}} a= a \sqrt{\frac{F_N}{N}},\ee
which we call nonlinear unitary transformation,
where $[a,a^+]=1.$

\subsubsection{Golden Angular Momentum}
Double Golden Oscillator algebra $su_F (2)$, determines the Golden Quantum angular momentum operators, defined as
$$J_+^F= b_1^+ b_2,\,\,\,\,\,\,J_-^F= b_2^+ b_1,\,\,\,\,\,\,J_z^F=\frac{N_1-N_2}{2},$$
 and satisfying commutation relations
\be [J_+^F,J_-^F]=(-1)^{N_2} F_{2J_z}= -(-1)^{N_1} F_{-2J_z},\ee
\be [J_z^F,J_{\pm}^F]=\pm J_{\pm}^F,\ee
where the Binet-Fibonacci operator is $$F_{N}= \frac{\varphi^N - (-\frac{1}{\varphi})^N}{\varphi+ \frac{1}{\varphi}}= [N]_F.$$
The Golden  quantum angular momentum operators $J_{\pm}^F$ may be written in terms of Fibonacci sequence and standard quantum angular momentum operators $J_{\pm}$ as
\be J_+^F=J_+ \sqrt{\frac{F_{N_1+1}}{N_1+1}} \sqrt{\frac{F_{N_2}}{N_2}}= \sqrt{\frac{F_{N_1}}{N_1}}\sqrt{\frac{F_{N_2+1}}{N_2+1}}J_+\ee
\be J_-^F=J_- \sqrt{\frac{F_{N_1}}{N_1}} \sqrt{\frac{F_{N_2+1}}{N_2+1}}=\sqrt{\frac{F_{N_1+1}}{N_1+1}}\sqrt{\frac{F_{N_2}}{N_2}}J_- .\ee
The Casimir operator for Binet-Fibonacci case is
\bea C^F&=& (-1)^{-J_z} \left(F_{J_z} F_{J_z+1} +(-1)^{-N_2} J_-^F J_+^F\right) \nonumber \\
&=& (-1)^{-J_z} \left(-F_{J_z}F_{J_z-1} +(-1)^{-N_2} J_+^F J_-^F \right) .\eea
The angular momentum
operators $J_{\pm}^F$ and $J_z^F$ act on state $|j,m\ra_F:$
\be J_+^F |j,m\ra_F =\sqrt {F_{j-m} F_{j+m+1}}|j,m+1\ra_F,\label{golden1}\ee
\be J_-^F |j,m\ra_F =\sqrt{F_{j+m} F_{j-m+1}}|j,m-1\ra_F ,\label{golden2} \ee
\be J_z^F |j,m\ra_F =m |j,m\ra_F. \label{golden3} \ee
The eigenvalues of Casimir operator $C_j^F$ are determined by product of two successive Fibonacci numbers
$$C_j^F= (-1)^{-j} F_{j} F_{j+1},$$
then the asymptotic ratio of two successive eigenvalues of Casimir operator gives Golden Ratio

$$\lim_{j \rightarrow \infty} \frac{(-1)^{-j} F_j F_{j+1}}{(-1)^{-j+1} F_{j-1} F_j}= - \varphi^2. $$

We can also construct representation  of our $F$-deformed angular momentum algebra in terms of double Golden boson
representation $b_1,b_2.$ The actions of $F$-deformed angular momentum operators to the state $|n_1,n_2\ra_F$ are given as follows :

\be J_+^F |n_1,n_2 \ra_F = b_1^{+} \, b_2 |n_1,n_2 \ra_F=
\sqrt{F_{n_1+1} F_{n_2}} |n_1+1,n_2-1 \ra_F,\ee
\be J_-^F |n_1,n_2 \ra_F=b_2^{+} \, b_1 |n_1,n_2 \ra_F=
\sqrt{F_{n_1} F_{n_2+1}}  |n_1-1,n_2+1 \ra_F , \ee
\be J_z^F |n_1,n_2 \ra_F =\frac{1}{2} (N_1-N_2)|n_1,n_2 \ra_F=
\frac{1}{2}(n_1-n_2)  |n_1,n_2 \ra_F. \ee
 The above expressions reduce to the familiar ones (\ref{golden1})-(\ref{golden3}) provided we define
 $$j\equiv \frac{n_1+n_2}{2},\,\,\,\,\,\,\,\,\, m\equiv \frac{n_1-n_2}{2}$$ $$|n_1,n_2\ra_F \equiv |j,m\ra_F,$$ and substitute
 $$n_1\rightarrow j+m, \,\,\,\,\,\,\,\,\,\,\,\,n_2\rightarrow j-m.$$

\subsubsection{Symmetrical $su_{i \varphi} (2)$ Quantum Algebra}
As an example of symmetrical $q$-deformed $su_q(2)$ algebra we choose the base as $q_i=i \varphi$ and $q_j=i \frac{1}{\varphi},$ then
our complex equation for base becomes
$$(i \varphi)^2 = i (i \varphi)-1.$$ The $\varphi$- deformed symmetrical angular momentum  operators remain the same as $J_{\pm}^{(s)}, J_z^{(s)}.$ The symmetrical quantum algebra with base $(i \varphi, \frac{i}{\varphi})$ becomes
\be [J_+^{\varphi},J_-^{\varphi}]= [2 J_z]_{\frac{i}{\varphi}}= [2 J_z]_{i \varphi,\frac{i}{\varphi}} (-1)^{(\frac{1}{2}-J_z)}, \ee
where $$[2 J_z]_{\frac{i}{\varphi}}= \frac{\varphi^{2 J_z}-\varphi^{-2 J_z}}{\varphi-\varphi^{-1}}$$ and
\be [J_z^{(s)}, J_{\pm}^{(s)}]= \pm J_{\pm}^{(s)}.\ee
\subsubsection{$\tilde{su}_F(2)$ Algebra}
One of the special cases of symmetrical $\tilde{su}_{(q,Q)}(2)$ algebra is constructed by choosing Binet-Fibonacci case $(q_i=\varphi,\,\,\, q_j=-\frac{1}{\varphi}).$ The generators of $\tilde{su}_F(2)$ algebra $\tilde{J}_{\pm}^{\varphi}, \tilde{J}_z^{\varphi}$
are given in terms of double bosons $b_1,\,b_2$ as follows :
\be \tilde{J}_+^F = (-1)^{-\frac{N_2}{2}} b_1^+ b_2, \ee
\be \tilde{J}_-^F  = b_2^+ b_1 (-1)^{-\frac{N_2}{2}}, \ee
\be \tilde{J}_z^F = J_z.\ee
satisfying  anti-commutation relation
\be  \tilde{J}_+^F \tilde{J}_-^F+ \tilde{J}_-^F \tilde{J}_+^F=\{\tilde{J}_+^F ,\tilde{J}_-^F\} = [2 J_z]_F,\ee
and $[\tilde{J}_z^F,\tilde{J}_{\pm}^F]= \pm \tilde{J}_{\pm}^F.$
The Casimir operator is written in the following forms
\bea \tilde{C}^F &=& (-1)^{J_z} \{F_{j_z} F_{j_z +1}-\tilde{J}_-^F \tilde{J}_+^F\}\nonumber \\
&=& (-1)^{J_z} \{\tilde{J}_+^F \tilde{J}_-^F- F_{j_z}F_{j_z -1} \}.\eea
The actions of the $F$-deformed angular momentum operators to the states $|j,m\ra_F $ are
\be \tilde{J}_+^F|j,m\ra_F= (-1)^{\frac{j-m}{2}} \sqrt{F_{j-m} F_{j+m+1}} |j,m+1\ra_F,\ee
\be \tilde{J}_-^F |j,m\ra_F= (-1)^{\frac{j-m}{2}} \sqrt{F_{j+m} F_{j-m+1}} |j,m-1\ra_F,\ee
\be \tilde{J}_z^F |j,m\ra_F= m |j,m\ra_F.\ee
And the  eigenvalues of Casimir operators are given by
\bea \tilde{C}^F |j,m\ra_F &=& \{(-1)^m F_m F_{m+1}- (-1)^j F_{j-m} F_{j+m+1}\}|j,m\ra_F \nonumber \\
&=& \{(-1)^j F_{j-m+1} F_{j+m} -(-1)^m F_m F_{m-1}\} |j,m\ra_F .\nonumber \eea

\end{document}